\documentclass[a4paper,12pt]{article}
\usepackage{amsmath,amsthm,amssymb,latexsym}
\title{{\bf \boldmath$L^p$ and weak \boldmath$L^1$ estimates for the maximal Riesz transform and the maximal Beurling transform}}
\author{\Large{\Large Joan Mateu and Joan Verdera}}
\newtheorem{teo}{Theorem}

\newtheorem{lemma}[teo]{Lemma}

\theoremstyle{definition}
\newtheorem*{gracies}{Acknowledgements}

\newcommand{\Rn}{{\mathbb R}^n}

\newcommand{\C}{\mathbb{C}}

\begin{document}
\date{10/02/2006}
\maketitle

\begin{abstract}
We prove $L^p$ estimates for the maximal Riesz transform in terms of the Riesz
transform itself, for $1<p \leq \infty$. We show that the corresponding weak~$L^1$~ estimate fails for the maximal Riesz transform, but surprisingly does hold for the
maximal Beurling transform.
\end{abstract}

\section{Introduction}

Let $T$ be a classical Calder\'{o}n-Zygmund singular integral operator
on $\Rn$ with smooth kernel $K$ and let $T^{\star}$ be the
associated maximal singular integral
$$
T^{\star}f(x)= \sup_{\epsilon > 0} | T^{\epsilon}f(x)| \,,
$$
where $T^{\epsilon}f(x)$  is the truncation at level $\epsilon$
$$
T^{\epsilon}f(x)= \int_{| y-x| > \epsilon} f(x-y) K(y) \,dy\, .
$$

In this note we consider the problem of controlling $T^{\star} f$
by $Tf$. The control we seek is in terms of $L^p$ and weak $L^1$
estimates and we restrict our attention to the Riesz and Beurling
transforms. In the forthcoming paper \cite{MPV} one considers
pointwise estimates and more general kernels.

The $j$-th Riesz transform, $1\leq j \leq n$, is the singular integral operator
\begin{equation}\label{eq1}
R_jf(x)= PV \int f(x-y) \frac{y_j}{| y | ^{n+1}} \,dy \equiv
\lim_{\epsilon\rightarrow 0} \int_{|x-y|> \epsilon} f(x-y)
\frac{y_j}{| y | ^{n+1}}\, dy \,.
\end{equation}

The principal value integral above exists for all $x$ if $f$ is a compactly supported
smooth function and one shows for such functions the $L^p$ estimate
$$
\| R_jf  \|_p \leq  C \| f \|_p , \quad 1<p<\infty\,,
$$
for some positive constant $C$ independent of $f$. Then
a bounded operator $R_j$ can be defined on $L^p$ by the obvious
density argument. A subtle issue arises when one tries to show
that the principal values in \eqref{eq1} exist for almost all $x$ in $\Rn$
for each function $f$ in $L^p$. Following a well known principle,
one looks for an $L^p$ estimate for the maximal Riesz transform
$R^{\star}_j f$ and one indeed proves that
$$
\| R^{\star}_jf  \|_p \leq  C \| f \|_p ,\quad 1<p<\infty\,.
$$
For the classical results just mentioned one may consult \cite{Du} or \cite{St}.
Our first result improves the last inequality.

\begin{teo} \label{T1}
For $1<p \leq\infty$ there exists a constant $C=C_{p,n}$ such that
\begin{equation}\label{eq2}
\| R^{\star}_j f  \|_p \leq  C \| R_j f \|_p\,,
\end{equation}
for each function $f$ belonging to some $L^q(\Rn)$, $1\leq q \leq \infty$.
\end{teo}

One may wonder whether \eqref{eq2} still holds in the limiting case $p=1$ and the
answer is provided by our next result.

\begin{teo} \label{T2}
Given $j$, $1 \leq j \leq n$, and a positive constant $C$, there exist a function~$f$
in $L^1(\Rn)$ such that
$$
\| R^\star_j f  \|_{1,\infty} \geq  C \| R_j f \|_1\,.
$$
\end{teo}

Notice that the above statement is not obvious even for the Hilbert transform, which
is the case $n=1$. The surprising fact is that the weak $L^1$ inequality does hold for
the Beurling transform.

The Beurling transform is the singular integral in the plane $\C$ defined by
$$
Bf(z)= \frac{1}{\pi}\, PV \int f(z-w)\frac{1}{w^2}\, dA(w)\,,
$$
where $dA$ denotes $2$-dimensional Lebesgue measure. The Beurling transform is
invertible on $L^p(\C)$, $1<p<\infty$, and the inverse is the singular integral with
kernel~$ 1 / \pi\,\overline{z}^2 $. Thus $\| f \|_p  \leq C\| Bf \|_p$, and
consequently, $\|B^\star f \|_p \leq C\| Bf \|_p$. Hence inequality~\eqref{eq2} with
$R_j$ replaced by $B$ is trivially true for $1< p< \infty$. Incidentally, let's remark
that the same can be said about \eqref{eq2} for $n=1$, because the Hilbert transform
is invertible on $L^p$, $1<p<\infty$.

\begin{teo}\label{T3}
We have
$$
\|B^\star f \|_{1,\infty} \leq C\| Bf \|_1 \,,
$$
for any $f$ belonging to some $L^q(\C)$, $1 \leq\ q \leq \infty$.
\end{teo}

Therefore, there is an astonishing difference between even and odd kernels. This
phenomenon will be further studied in \cite{MPV}. One should not think that
Theorem~\ref{T3} is due to an extremely special property of the Beurling transform.
 Indeed, it can be shown that it holds for a significant class of even kernels,
 although it is not yet
clear to the authors what are the best possible assumptions one should require.

Taking into account the classical estimate $\|B^\star f \|_{1,\infty} \leq C\,\| f
\|_1 \,,$ Theorem $3$ entails
$$\|B^\star f \|_{1,\infty} \leq C\, {\rm min}\{\| Bf \|_1, \|f \|_1 \}\, ,\, f\in L^1(\C)\,. $$
\noindent Notice that there are functions $f\in L^1(\C)$, $\|f\|_1 =1$, such that
$\|Bf \|_1$ is as small as desired. On the other hand, one can also find a function $f
\in L^1(\C)$ with $\|Bf \|_1 = 1$ and $\|f\|_1$ very small. Thus neither the
inequality in Theorem $3$ nor the classical weak $L^1$ estimate for $B^\star f$ is
stronger than the other.

 In the next three sections we will provide the proofs of the above three theorems.
We adhere to the standard convention of denoting by $C$ a positive constant,
independent of the relevant parameters involved, and which may vary from an occurrence
to another. Our notation and terminology are standard. For example, $A \simeq B$ means
that the two quantities $A$ and $B$ satisfy the relation $C^{-1} A \leq B \leq C A$,
for some constant $C\geq1$.

The problem considered in this paper arose when the second named author was working on
the question, still open,  of whether the Riesz kernels characterize uniform
rectifiability in dimensions greater than $1$ (see \cite[p.~139]{DS}; see also
\cite{Ve} for a survey about the one dimensional problem and related results).

 The authors are grateful to Carlos P\'{e}rez for an illuminating conversation
 which lead to Theorem~\ref{T2}.

\section{Proof of Theorem 1}

In this section we may assume that $n > 1$, because the result is obviously true in
dimension~$1$. Let $M$ be the Hardy-Littlewood maximal operator. It is enough to prove
that for each $s>1$ there exists a positive constant $C_s$ such that
$$
R_j^\star f (x) \leq C_s  {M({| R_jf|}^s)}^\frac{1}{s} (x), \quad x\in \Rn\,.
$$

Let $B$ denote the unit ball centered at the origin, $B^c$
its complement in $\Rn$ and let $K_j(x)= \frac{x_j}{|x|^{n+1}}$ be the kernel of the
$j$-th Riesz transform.

\begin{lemma}\label{lem4}
There exists a function $h$ such that
\begin{equation}\label{eq3}
\chi_{B^c}(x)\,K_j(x)= R_jh(x),\quad x \in \Rn, \quad 1\leq j \leq n\, .
\end{equation}
\end{lemma}

\begin{proof}
As it is well known
$$
\partial_j \left(\frac{1}{|x|^{n-1}}\right)= -(n-1)K_j(x)\,,
$$
in the distributions sense. Consider the function $\varphi$ that takes the
value $1$ on $B$ and $1/|x|^{n-1}$ on $B^c$ . Since $\varphi$ is continuous on the
boundary of $B$, we have, in the distributions sense,
$$
\partial_j \varphi = -(n-1)\chi_{B^c} K_j \,.
$$

Taking the Fourier transform one sees that for an appropriate
constant $c_n$ one has
$$
\varphi= \frac{1}{|x|^{n-1}}\star c_n \sum_{i=1}^n R_i(\partial_i \varphi)  \,.
$$
Therefore \eqref{eq3} follows by taking the $j$-th derivative of
$\varphi$ and setting, for another suitable constant $c_n$,
\begin{equation}\label{eq4}
 h= c_n \sum_{i=1}^n R_i(\chi_{B^c}K_i)\,.
\end{equation}
\end{proof}

Set
$$ R_j^\epsilon f(x)= \int_{|y|>\epsilon} f(x-y) \frac{y_j}{|y|^{n+1}}\, dy\,. $$
 We have to show that
$$ | R_j^\epsilon f(x)| \leq C_s \,
{M({| R_jf|}^s)}^\frac{1}{s} (x), \quad x\in \Rn\,.$$
By translation and dilation invariance we can assume, without loss of
generality, that $x=0$ and $\epsilon = 1$. Then
\begin{equation*}
\begin{split}
 R_j^1 f(0) &=  - \int  \chi_{B^c}(x)\,K_j(x) \,f(x)\, dx = - \int R_jh(x)\,f(x)\,dx\\*[5pt]
& = \int h(x) \, R_jf(x)\, dx\,.
\end{split}
\end{equation*}
We will see below that $h \in L^1(\Rn)$, but $h$ has not an integrable
decreasing radial majorant. To overcome this difficulty we split the last integral
into two pieces
$$
R_j^1 f(0) = \int_{2B} h(x) \, R_jf(x) \,dx + \int_{(2B)^c} h(x) \, R_jf(x) \,dx \equiv
I_1+I_2 \,.
$$

To estimate the term $I_1$ notice that $h$ belongs to $L^q(2B)$ for $1 \leq q <
\infty$. This follows from~\eqref{eq4}, because the functions $\chi_{B^c}\,K_i$ are in
$L^q(\Rn)$ for $1<q$, and the $R_i$ are bounded on $L^q(\Rn)$ if $1< q<\infty$. Thus,
if $s'$ denotes the exponent conjugate to $s$, by Holder's inequality we obtain
$$
|I_1|\leq \left(\int_{2B} |R_jf |^s\, dx\right)^{\frac{1}{s}} \left(\int_{2B} | h |^{s'} \,dx
\right)^{\frac{1}{s'}} \leq C_s \; M(|R_jf|^s)^\frac{1}{s} (0)\,.
$$

The term $I_2$ can easily be estimated if we first prove that
\begin{equation}\label{eq5}
|h(x)| \leq C \frac{1}{|x|^{n+1}},\quad |x|\geq 2\,.
\end{equation}
Indeed, the preceding decay inequality yields
$$
|I_2| \leq  C  \int_{(2B)^c} |R_jf (x)| \frac{1}{|x|^{n+1}}\,
dx \leq  C \, M(R_jf)(0)\,,
$$
which is not greater than $ C \, M(|R_jf|^s)^\frac{1}{s} (0)$ by Holder's
inequality.

To prove \eqref{eq5} express $h/c_n$ as
$$
\frac{h}{c_n}= \sum_{i=1}^n K_i \star \chi_{B^c} K_i = \sum_{i=1}^n K_i \star K_i -
\sum_{i=1}^n K_i \star \chi_{B} K_i = c'_n  \delta_0 -  \sum_{i=1}^n R_i
(\chi_{B} K_i)\,.
$$
If $|x|> 1$ we have
\begin{equation*}
\begin{split}
 R_i (\chi_{B} K_i)(x) &=  \lim_{\epsilon \rightarrow 0} \int_{\epsilon< |y|< 1}
K_i(x-y) K_i(y)\, dy \\*[5pt]
&=  \lim_{\epsilon \rightarrow 0}
\int_{\epsilon< |y|< 1} (K_i(x-y)- K_i(x)) K_i(y)\, dy \,.
\end{split}
\end{equation*}
Since
$$
| K_i(x-y)- K_i(x)| \leq  C  \frac{|y|}{|x|^{n+1}}, \quad |x|\geq
2,\quad  |y|\leq 1 \,,
$$
we obtain, for $|x|\geq 2$,
$$
|R_i (\chi_{B} K_i)(x)|\leq  C  \int_{|y|<1}
\frac{1}{|x|^{n+1}}  \frac{1}{|y|^{n-1}} \,dy =
\frac{C}{|x|^{n+1}}\,,
$$
which gives \eqref{eq5} and completes the proof of Theorem~\ref{T1}.

\section{Proof of Theorem 2}
We prove Theorem~\ref{T2} for $j=1$. Set $b=(1,0,\dotsc,0)$ and
$a=(-1,0,\dotsc,0)$ and let~$\mu$~be the length measure on the
segment joining  $a$ and $b$. For an appropriate constant $c_n$ we
have
$$
\mu = c_n  \left(\frac{1}{|x|^{n-1}} \star \sum_{j=1}^n R_j(\partial_j \mu)\right)\,,
$$
as one can easily see by taking Fourier transforms on both sides. For $n=1$
one should replace $\frac{1}{|x|^{n-1}}$ by $\log |x|$ in the formula above.
We have
$$
\delta_a - \delta _b = \partial_1 \mu = R_1 \left(c_n \sum_{j=1}^n
R_j(\partial_j \mu)\right)\,.
$$
Set
$$
T = c_n \sum_{j=1}^n R_j(\partial_j \mu)\,,
$$
so that
\begin{equation}\label{eq6}
\delta_a - \delta _b = R_1(T)\,.
\end{equation}

Let $\varphi$ be a non-negative continuously differentiable function with compact
support contained in the unit ball $B$ such that $\int \varphi = 1$, and set
$$
\varphi_\epsilon(x) = \frac{1}{\epsilon^n} \varphi \left(\frac{x}{\epsilon}\right)\,.
$$
Convolving the identity~\eqref{eq6} with $\varphi_\epsilon$ we obtain
$$
\varphi_\epsilon(x-a) - \varphi_\epsilon(x-b) = R_1(T \star
\varphi_\epsilon)\,.
$$
Let $f_\epsilon$ stand for $ T \star \varphi_\epsilon $, so that
$$
\| R_1(f_\epsilon)\|_1 \leq 2 \,.
$$
Now
$$
f_\epsilon =  c_n \sum_{j=1}^n R_j(\mu \star \partial_j \varphi_\epsilon)
$$
and $\mu \star \partial_j \varphi_\epsilon$ is a compactly supported function in
$L^\infty(\Rn)$ with zero integral. Thus $\mu \star \partial_j \varphi_\epsilon$ is a
function in the Hardy space $H^1(\Rn)$ (in fact, \,it is a multiple of an atom) and so
$f_\epsilon \in L^1(\Rn)$.  Hence we only need to show that
\begin{equation}\label{eq7}
\| R_1^{\star} f_\epsilon \|_{1,\infty} \geq  \frac{1}{C}
\log\left(\frac{1}{\epsilon}\right) \,,
\end{equation}
for some constant $C\geq 1$ and for sufficiently small
positive $\epsilon$.

Recall from Lemma~\ref{lem4} that there is a function $h$ such that $R_1(h)= \chi_{B^c}
 K_1$. Dilating we get $R_1(h_\delta)= \chi_{(\delta B)^c} K_1 $,
 where $h_\delta(x)= \frac{1}{\delta^n} h(\frac{x}{\delta})$ and
 $\delta B$ is the ball of radius $\delta$ centered at the origin.
 Therefore
\begin{equation*}
\begin{split}
R_1^\delta(T)(x)&=  (R_1(h_\delta) \star T) (x) = - (h_\delta
\star R_1(T))(x)\\*[5pt]
&=  h_\delta \star (\delta_b - \delta_a)(x)=
h_\delta(x-b)-h_\delta(x-a)
\end{split}
\end{equation*}
and, convolving with $\varphi_\epsilon$,
\begin{equation}\label{eq8}
 R_1^\delta(f_\epsilon)(x)= \left(\left(h_\delta(x-b)-h_\delta(x-a)\right) \star
\varphi_\epsilon \right) (x)\,.
\end{equation}
To go further we need to understand the singularity of
$h$.

\begin{lemma}\label{lem5}
We have, for some constant $c_0$,
$$
h(x)= b(x)+ c_0\, p(x)
$$
where $|b(x)| \leq C$, $x \in \Rn$, and
$$
p(x)= \int_{|y|=1} \frac{1}{|x-y|^{n-1}} \,d\sigma (y)\,,
$$
$\sigma$ being the $n-1$ dimensional surface measure on the unit sphere.
\end{lemma}

\begin{proof}
Assume first that $|x|<1$. Let $\omega_j = (-1)^{j-1}\, dx_1
\wedge\dotsb\wedge dx_{j-1}\wedge dx_{j+1}\wedge\dotsb \wedge dx_n $.
Apply Green-Stokes to the domain $ 1< |x|<R$ and then let $R
\rightarrow \infty$ to obtain
\begin{equation}\label{eq9}
\begin{split}
& \sum_{j=1}^{n}  \int_{|y|= 1}
\frac{1}{|x-y|^{n-1}}\partial_j\left(\frac{1}{|y|^{n-1}}\right)\omega_j  \\*[5pt]
=c_1 &\sum_{j=1}^{n}  R_j(\chi_{B^c}\,K_j)(x) + c_2 \int_{|y|>1}
\frac{1}{|x-y|^{n-1}} \frac{1}{|y|^{n+1}} \,dy \,,
\end{split}
\end{equation}
for some constants $c_1$ and $c_2$. Thus, with other
constants $c_3$ and $c_4$,
$$
h(x)= c_3   \int_{|y|= 1} \frac{1}{|x-y|^{n-1}} \sum_{j=1}^{n}
y_j \,\omega_j + c_4 \int_{|y|>1} \frac{1}{|x-y|^{n-1}}
\frac{1}{|y|^{n+1}} \,dy \,.
$$
Now, the form $\sum_{j=1}^{n} y_j \,\omega_j $ is
invariant by rotations and so, understood as a measure, is a
constant multiple of $d\sigma$. To complete the proof of the lemma
we only need to show that the second term in the right hand side
of the preceding identity is a bounded function of $x$. This is
very easy if we split the domain of integration into two parts
according to whether $|x-y|<1$ or $|x-y|>1$. We then get
$$
\int_{|y|>1} \frac{1}{|x-y|^{n-1}} \frac{1}{|y|^{n+1}} \,dy \leq
\int_{|x-y|<1} \frac{1}{|x-y|^{n-1}} \,dy + \int_{|y|>1}
\frac{1}{|y|^{n+1}} \,dy \leq C \,.
$$

If $|x|>1$ the argument is basically the same except that one has to delete a small
ball centered at $x$ before applying Green-Stokes and then let the radius tend to
zero. We again get the identity~\eqref{eq9}, where now the Riesz transforms really involve
principal values.
\end{proof}

The next simple lemma describes precisely the singularity of the potential $p(x)$.
Denote by $d(x)= ||x|-1|$ the distance from $x$ to the unit sphere $\{|x|=1 \}$.

\begin{lemma}\label{lem6}
We have
$$
\int_{|y|=1} \frac{d\sigma(x)}{|x-y|^{n-1}}  \simeq \log
\frac{1}{d(x)}\,,
$$
provided $d(x) \leq \frac{1}{2}$.
\end{lemma}

\begin{proof}
Take $x$ such that $d(x)\leq \frac{1}{2}$ and set
$$
A_k = \{ y: |y|= 1 \quad \text{and}\quad  d(x)2^k \leq |y-x| \leq
d(x)2^{k+1} \},\quad 0\leq k\,.
$$
If $N$ is chosen appropriately, then
$$
\int_{|y|=1} \frac{d\sigma(y)}{|x-y|^{n-1}} = \sum_{k=0}^N \int_{A_k}
\frac{d\sigma(y)}{|x-y|^{n-1}} \simeq \sum_{k=0}^N \frac{\sigma(A_k)}{(d(x)2^k)^{n-1}}
\simeq \sum_{k=0}^N 1 \simeq N \simeq \log \frac{1}{d(x)}\,.
$$
\end{proof}

\begin{lemma}\label{lem7}
There exists a constant $C\geq1$ such that, for sufficiently small
$\epsilon>0$ we have
\begin{alignat*}{2}
(p_\delta \star \varphi_\epsilon)(x) &\geq \frac{1}{C}
\frac{1}{\delta^{n}}\log \frac{\delta}{\epsilon}\,, &\quad&\text{if}\quad \operatorname{dist}(x,\partial B(0,\delta)) < 2\epsilon\, \delta\,,\\*[5pt]
(p_\delta \star \varphi_\epsilon)(x) &\leq \frac{1}{C}
\frac{1}{\delta^{n}} \log \frac{\delta}{\operatorname{dist}(x,\partial
B(0,\delta))}\,, &\quad &\text{if}\quad 2\epsilon\,\delta < \operatorname{dist}(x,\partial B(0,\delta))  \leq \frac{\delta}{2} \,.
\end{alignat*}
\end{lemma}

\begin{proof}
Assume first that $\delta=1$ and that $d(x)= \operatorname{dist}(x,\partial
B(0,1))< 2\epsilon$\,. For each $y$ in the ball $B(x,\epsilon)$
one has $d(y)\leq d(x) +\epsilon \leq 3\epsilon$, and so
$$
\int \left(\log \frac{1}{d(y)}\right)\varphi_\epsilon(x-y)\, dy \geq \log
\frac{1}{ 3 \epsilon}\,,
$$
which is not greater than $\frac{1}{2} \log \frac{1}{\epsilon}$ if
$\epsilon$ is small enough.

If $2\epsilon < d(x) < \frac{1}{2}$, then
$$
\int \left(\log \frac{1}{d(y)}\right)\varphi_\epsilon(x-y)\, dy \leq \log
\frac{2}{ d(x)} \leq 2 \log \frac{1}{d(x)}\,.
$$

Consider now an arbitrary positive $\delta$. Then
$$
(p_\delta \star \varphi_\epsilon)(x)= \frac{1}{\delta^n}(p \star
\varphi_{\frac{\epsilon}{ \delta}} )\left(\frac{x}{\delta}\right)\,.
$$
Since
$$
\operatorname{dist}\left(\frac{x}{\delta},\partial B(0,1)\right)=
\operatorname{dist}(x,\partial B(0,\delta))\delta^{-1}\,,
$$
the lemma follows.
\end{proof}

We proceed now to prove \eqref{eq7}. Consider the cone $K$ with vertex at $b$, aperture
$\frac{\pi}{4}$ and with axis the positive $x_1$-axis. In other words,
$$
K= \left\{x \in \Rn:\langle x-b,b\rangle  \geq \frac{1}{\sqrt{2}}\;|x-b|\right\}\,.
$$
Take $x \in K$ and set $\delta = |x-b|$. We are going to apply the second
inequality in Lemma~\ref{lem7} with $x$ replaced by $x-a$. Thus we have to check that
\begin{equation}\label{eq10}
2\epsilon \delta < |x-a|-\delta \leq \frac{\delta}{2}\,.
\end{equation}
The second inequality is obvious if we assume $\delta \geq 4 $, because then
$|x-a|\leq |x-b|+2 \leq \frac{3}{2}\delta$. A simple estimate based on the fact that
$x \in K$  shows that $ |x-a|-\delta \geq \sqrt{2}$ and thus the first inequality
in~\eqref{eq10} holds provided that $\delta < \frac{1}{\sqrt{2} \epsilon}$. Using
\eqref{eq8}, Lemma~\ref{lem6} and Lemma~\ref{lem7}, we obtain
\begin{equation*}
\begin{split}
 |R_1^\delta(f_\epsilon)(x)|&= |((h_\delta(x-b)-h_\delta(x-a)) \star
\varphi_\epsilon) (x)| \\*[5pt]
&\geq \frac{1}{C}
|((p_\delta(x-b)-p_\delta(x-a)) \star \varphi_\epsilon) (x)| - C\\*[5pt]
&\geq \frac{1}{\delta^n} \left(\frac{1}{C} \log \frac{\delta}{\epsilon} - C
\log \frac {\delta}{\operatorname{dist}(x,\partial B(a,\delta))}- C\right)\,.
\end{split}
\end{equation*}
Since $\operatorname{dist}(x,\partial B(a,\delta))= |x-a|-\delta > \sqrt{2}>
1$,
$$
|R_1^\delta(f_\epsilon)(x)| \geq \frac{1}{\delta^n} \left(\frac{1}{C}
\log \frac{\delta}{\epsilon} - C \log \delta - C\right)\,,
$$
which is greater than or equal to
$$
\frac{1}{\delta^n}
\left(\frac{1}{C} \log \frac{1}{\epsilon} - C \log \delta\right)\,,
$$
because $\log\delta \geq \log 4 > 1$. Since
$$
C \log \delta
\leq\frac{1}{2C} \log\frac{1}{\epsilon}\,,
$$
provided
$$
\delta \leq\epsilon^{-\eta},\quad \eta = \frac{1}{2C^2}\,,
$$
we conclude that
\begin{equation}\label{eq11}
|R_1^\delta(f_\epsilon)(x)| \geq  \frac{1}{C} \frac{1}{|x-b|^n} \log
\frac{1}{\epsilon}, \quad x \in K,\quad 4\leq |x-b|=\delta \leq \epsilon^{-\eta}\,,
\end{equation}
for $\epsilon$ so small that $\epsilon^{-\eta} \leq \frac{1}{\sqrt{2} \epsilon}\,.$

Denote by $|E|$ the Lebesgue measure of the set $E$. By~\eqref{eq11}, for small $\epsilon$ we
 obtain
\begin{equation*}
\begin{split}
|\{x \in \Rn : R_1^{\star}f_\epsilon (x) \!>\! 1 \}| &\geq  \left|\left\{x \!\in\! K
: 4\leq |x-b| \leq \epsilon^{-\eta}\text{ and } \frac{1}{C}
\frac{1}{|x-b|^{n}} \log \frac{1}{\epsilon} > 1 \right\}\right| \\*[9pt]
&=  \left|\left\{x \!\in\! K : 4\leq |x-b| \leq \epsilon^{-\eta}\text{ and } |x-b| \!<\! \left(\frac{1}{C}\log \frac{1}{\epsilon}\right)^{\!\!\frac{1}{n}\!}\right\}\right| .
\end{split}
\end{equation*}
If $\epsilon$ is such that
$$
\left(\frac{1}{C} \log
\frac{1}{\epsilon}\right)^{\frac{1}{n}} \leq \epsilon^{-\eta}\,,
$$
then we get
$$
|\{x \in \Rn : R_1^{\star}f_\epsilon (x) > 1 \}| \geq \left|\left\{x \in K :
4\leq |x-b| \leq  \left(\frac{1}{C} \log
\frac{1}{\epsilon}\right)^{\frac{1}{n}}\right\}\right|\,.
$$
Taking $\epsilon$ small enough we can further assume that  $8 \leq
(\frac{1}{C} \log \frac{1}{\epsilon})^{\frac{1}{n}}$. Therefore
$$
|\{x \in \Rn : R_1^{\star}f_\epsilon (x) > 1 \}| \geq \frac{1}{C}
\log \frac{1}{\epsilon}\,,
$$
which yields
$$
\|R_1^{\star}f_\epsilon  \|_{1,\infty} \geq \frac{1}{C} \log
\frac{1}{\epsilon}\,,
$$
and completes the proof of Theorem~\ref{T2}.

\section{Proof of Theorem 3}

Theorem~\ref{T3} follows from the pointwise inequality
\begin{equation}\label{eq12}
B^{\star}f (z)\leq M(Bf)(z),\quad z \in \C\,.
\end{equation}

Let's remark, incidentally, that the preceding estimate is an improvement of Cotlar's
inequality for the Beurling transform
$$
B^{\star}f (z)\leq  C (M(Bf)(z)+ M(f)(z)),\quad z \in \C\,,
$$
 because the term $ Mf(z)$ does not appear in the right hand side of~\eqref{eq12}.

To prove \eqref{eq12} we first show the formula
\begin{equation}\label{eq13}
\frac{1}{z^2} \chi_{D^c}(z) = B(\chi_{D})(z)\,,
\end{equation}
where $D$ is the disc of center $0$ and radius $1$. Let
$$
C(f)(z)= \frac{1}{\pi} \int f(z-w) \frac{1}{w} \,dA(w)
$$
be the Cauchy transform of the function $f$, so that
$$
\frac{\partial}{\partial \overline{z}} C(f) = f \quad \text{and}\quad \frac{\partial}{\partial z} C(f) = - B(f)\,,
$$
in the distributions sense.

Consider the function $F(z)$ which takes the value $\overline{z}$ on $D$ and the value
$\frac{1}{z}$ on $D^c$. Then
$$
\frac{\partial F}{\partial \overline{z}}  = \chi_D\quad \text{and}
\quad \frac{\partial F}{\partial z}  = - \frac{1}{z^2}
\chi_{D^c}\,.
$$
Notice that $F= C(\chi_D)$, because $F-C(\chi_D)$ is an entire function
vanishing at $\infty$. Thus \eqref{eq13} holds.

Hence
\begin{equation*}
\begin{split}
B^1(f)(0)&=   \frac{1}{\pi} \int \chi_{D^c}(z) \frac{1}{z^2}
f(z)\,dA(z)\\*[5pt]
&=  \frac{1}{\pi} \int B(\chi_D)(z) f(z)\,dA(z)\\*[5pt]
&=  \frac{1}{\pi} \int_D B(f)(z)\,dA(z)\,.
\end{split}
\end{equation*}
Dilating and translating one obtains
$$
B^\epsilon(f)(z)= \frac{1}{\pi \epsilon^2} \int_{D(z,\epsilon)}
B(f)(w)\,dA(w)\,,
$$
and, consequently, \eqref{eq12}.

\begin{gracies}
 The authors were partially supported by grants\newline  2005SGR00774
(Generalitat de Catalunya),  MTM2004-00519  and HF2004-0208.
\end{gracies}

\begin{tabular}{l}
Joan Mateu\\
Departament de Matem\`{a}tiques\\
Universitat Aut\`{o}noma de Barcelona\\
08193 Bellaterra, Barcelona, Catalonia\\
{\it E-mail:} {\tt mateu@mat.uab.es}\\ \\
Joan Verdera\\
Departament de Matem\`{a}tiques\\
Universitat Aut\`{o}noma de Barcelona\\
08193 Bellaterra, Barcelona, Catalonia\\
{\it E-mail:} {\tt jvm@mat.uab.es}
\end{tabular}

\end{document}